\documentclass[12pt]{amsart}
\usepackage[margin=1in]{geometry}
\usepackage[dvipsnames,svgnames,x11names]{xcolor}
\usepackage[hyphens]{url}
\usepackage[colorlinks=true,
            linkcolor=Maroon,
            citecolor=blue,
            urlcolor=blue,
            hypertexnames=false,
            linktocpage]{hyperref}
\usepackage{bookmark}
\usepackage{amsmath,thmtools,mathtools,amssymb}
\mathtoolsset{showonlyrefs=true}
\allowdisplaybreaks
\newcommand{\R}{\mathbb{R}}
\newcommand{\Sn}{S^{n-1}}
\newcommand{\one}{\mathbf{1}}

\newcommand{\ip}[2]{\langle #1,#2\rangle}
\newcommand{\norm}[1]{\left\|#1\right\|}
\newcommand{\Vol}{\operatorname{Vol}}
\newcommand{\relint}{\operatorname{relint}}
\newcommand{\hess}{\mathop{\overline{\operatorname{Hess}}}\nolimits}

\theoremstyle{plain}
\newtheorem{theorem}{Theorem}[section]
\newtheorem{proposition}[theorem]{Proposition}
\newtheorem{lemma}[theorem]{Lemma}
\newtheorem{corollary}[theorem]{Corollary}

\theoremstyle{definition}
\newtheorem{definition}[theorem]{Definition}

\theoremstyle{remark}
\newtheorem{remark}[theorem]{Remark}

\numberwithin{equation}{section}

\newcommand{\eq}[1]{\begin{equation}\begin{alignedat}{2}#1\end{alignedat}\end{equation}}

\makeatletter
\renewcommand{\subsection}{\@startsection{subsection}{2}%
  \z@{.5\linespacing\@plus.7\linespacing}{.25\linespacing}%
  {\normalfont\bfseries}}
\makeatother

\title[Continuity of the centro-affine spectral gap]
  {Continuity of the centro-affine\\ spectral gap for polytopes}
\author[Y. Hu, M. N. Ivaki]{Yingxiang Hu, Mohammad N. Ivaki}

\begin{document}
\begin{abstract}
For $n\geq 3$, we prove that if full-dimensional origin-symmetric polytopes $P_{\nu}$ in $\R^n$ converge in the Hausdorff metric to a full-dimensional polytope $P$, then their first nonconstant even centro-affine eigenvalues satisfy $\lambda_{1,e}(P_{\nu})\to \lambda_{1,e}(P)$.
\end{abstract}
\maketitle

\section{Introduction}

Throughout the paper, $n\geq 3$. The notation follows \cite{HI26}, where $\lambda_{1,e}(P)$ denotes the first nonconstant even eigenvalue of the weighted graph Laplacian on the facets of $P$.

\begin{theorem}\label{thm:polyhedral-lower-semicontinuity}
Suppose that full-dimensional origin-symmetric polytopes $P_{\nu}\subseteq\R^n$ converge in the Hausdorff metric to a full-dimensional polytope $P$. Then
\eq{
\lambda_{1,e}(P_{\nu})\to \lambda_{1,e}(P).
}
\end{theorem}

To prove this theorem, we construct a test function on $P$ from centered unit-mass even eigenfunctions $f_{\nu}$ for $\lambda_{1,e}(P_{\nu})$. Extending each $f_{\nu}$ radially defines functions $\hat{f}_{\nu}$ on the common domain $\R^n$. Along a subsequence realizing a finite $\liminf$ of $\lambda_{1,e}(P_{\nu})$, local $\mathrm{BV}$ estimates on $\hat{f}_{\nu}\lvert\hat{f}_{\nu}\rvert$ yield a further subsequence of $(\hat{f}_{\nu})$ converging strongly in $L^2_{\mathrm{loc}}(\R^n)$ (\autoref{lem:sec-3-radial-compactness}). A dual energy formula for radial extensions, analogous to the perimeter formula in \cite{Mag12}, shows that the $L^2_{\mathrm{loc}}$ limit is the radial extension $\hat{f}$ of a facet function $f$ on $P$ with energy at most $\liminf_{\nu\to \infty}\lambda_{1,e}(P_{\nu})$. Moreover, $f$ is even, with zero weighted mean and unit mass (\autoref{prop:sec-3-facet-function-compactness}). Evaluating the Rayleigh quotient of $f$ proves the lower semicontinuity claim in \autoref{thm:polyhedral-lower-semicontinuity}. The upper semicontinuity is a consequence of \cite[Lem. 8.9]{HI26}. See also \autoref{cor:radial-eigenfunction-convergence} for subsequential convergence of first nonconstant even eigenfunctions on $P_{\nu}$ to a first nonconstant even eigenfunction on $P$.

The dual energy formula also implies invariance of $\lambda_{1,e}$ under invertible linear maps (\autoref{lem:linear-invariance}). Together with John's theorem and Blaschke's selection theorem, \autoref{lem:linear-invariance} and the lower semicontinuity from \autoref{thm:polyhedral-lower-semicontinuity} remove the local minimality assumption in \cite[Thm. 1.1]{HI26} and establish $\lambda_{1,e}(P)\geq n-1+1/(n-1)$ for every full-dimensional origin-symmetric polytope $P\subseteq\R^n$ (\autoref{cor:polyhedral-spectral-gap}). Consequently, the local and global $L^{p_n}$-Brunn--Minkowski inequalities hold at $p_n=1-1/(n-1)$ for all origin-symmetric convex bodies in $\R^n$ (\autoref{thm:local-Lpn-inequality}, \autoref{cor:global-Lpn-inequality}), and every origin-symmetric $C^2_{+}$ convex body $K\subseteq\R^n$ satisfies $\lambda_{1,e}(K)\geq n-1+1/(n-1)$. We point out that the logarithmic Brunn--Minkowski inequality was recently proved in \cite[Thm. 1.1]{OAI26}. In particular, the sharp lower bound $\lambda_{1,e}(K)\geq n$ holds for every origin-symmetric convex body $K$ that is a polytope or belongs to $C^2_{+}$.

\section{Background}
For $x,y\in\R^n$, write $|x|$ for Euclidean length and $\ip{x}{y}$ for the Euclidean inner product. The notation $B_2^n$ stands for the closed Euclidean unit ball centered at the origin, $\Sn$ for the unit sphere, and $\mathcal{H}^k$ for $k$-dimensional Hausdorff measure on $\R^n$. We write $\Vol=\mathcal{H}^n$ for $n$-dimensional volume. 

A \emph{convex body} is a compact convex set in $\R^n$ with nonempty interior. A set $K\subseteq\R^n$ is \emph{origin-symmetric} if $K=-K$. 

For a nonempty set $E\subseteq\R^n$, we denote its affine hull by $\operatorname{aff}E$, its convex hull by $\operatorname{conv}E$, and its affine dimension by $\dim E=\dim\operatorname{aff}E$. Unless otherwise specified, the interior $\operatorname{int}E$ and boundary $\partial E$ are taken in $\R^n$, while the relative interior $\relint E$ and relative boundary $\operatorname{relbd}E$ are taken in $\operatorname{aff}E$. For sets $K,L\subseteq\R^n$ and $t\geq 0$, the Minkowski sum of $K$ and $tL$ is $K+tL=\{x+ty:x\in K,\ y\in L\}$.

The \emph{support function} of a nonempty compact set $K\subseteq\R^n$ is
\eq{
h_K(u)=\max_{x\in K}\ip{x}{u}
\quad(u\in\R^n).
}

A \emph{polytope} $P\subseteq\R^n$ is the convex hull of a finite nonempty subset of $\R^n$. As in \cite[Sec. 2]{HI26}, a face of $P$ is \emph{proper} if it is neither $\varnothing$ nor $P$. For a full-dimensional polytope, faces of dimension $0$ are vertices, of dimension $1$ are edges, of dimension $n-2$ are ridges, and of dimension $n-1$ are facets.

If a convex body $K$ contains the origin in its interior, its \emph{polar body} is
\eq{
K^{\ast}=\{y\in\R^n:\ip{x}{y}\leq 1\text{ for every }x\in K\}.
}
For a full-dimensional polytope $P$ containing the origin in its interior, polarity reverses proper face inclusions and sends facets and ridges of $P$ to vertices and edges of $P^{\ast}$, respectively; see \cite[p. 109, (2.28)]{Sch14}.

Suppose that $P\subseteq\R^n$ is full dimensional and contains the origin in its interior. A complete irredundant facet representation indexed by a finite set $J$ has the form
\eq{
P=\bigcap_{i\in J}\{x\in\R^n:\ip{x}{u_i}\leq h_i\},
\quad
|u_i|=1,
\quad
h_i> 0,
}
with one inequality for each facet of $P$. The vector $h=(h_i)_{i\in J}$ is the \emph{support vector}, and the facet with label $i$ is
\eq{
F_i=\{x\in P:\ip{x}{u_i}=h_i\}.
}

\begin{definition}
The \emph{mass} of the facet $F_i$ is defined as $m_i=h_i\mathcal{H}^{n-1}(F_i)$. Two facets are \emph{adjacent} when they meet in a ridge. For adjacent facets $F_i$ and $F_j$, define
\eq{
R_{ij}=F_i\cap F_j,
\quad
\theta_{ij}=\arccos\ip{u_i}{u_j},
\quad
c_{ij}=\frac{h_ih_j\mathcal{H}^{n-2}(R_{ij})}{\sin\theta_{ij}}.
}
For adjacent facets $F_i$ and $F_j$, we have $0<\theta_{ij}<\pi$. The number $c_{ij}$ is the \emph{ridge conductance}, with $c_{ij}=0$ for nonadjacent pairs $\{i,j\}$. The \emph{labeled facet network} has vertex set $J$, edges joining adjacent facets, vertex weights $m_i$, and edge weights $c_{ij}$.
\end{definition}

The facet-adjacency graph of $P$ is connected \cite[Lem. 2.2]{HI26}.

A \emph{facet function} is a vector $f=(f_i)_{i\in J}\in\R^J$. We write $\one$ for the facet function equal to one at every label. For facet functions $f,g\in\R^J$, the mass and energy forms are
\eq{
\mathsf{M}_P(f,g)&=\sum_{i\in J}m_if_ig_i,\\
\mathsf{L}_P(f,g)&=\sum_{\{i,j\}\in\binom{J}{2}}c_{ij}(f_i-f_j)(g_i-g_j).
}
The second sum runs over the set $\binom{J}{2}$ of two-element subsets of $J$, with each unordered pair counted once. A facet function is said to be \emph{centered} when $\mathsf{M}_P(f,\one)=0$. The \emph{weighted mean} of $f$ is defined by
\eq{
\bar{f}=\frac{\mathsf{M}_P(f,\one)}{\mathsf{M}_P(\one,\one)}=\frac{\sum_{i\in J}m_if_i}{\sum_{i\in J}m_i}.
}
Write
\eq{
\widetilde{\mathsf{M}}_P(f,g)
=\mathsf{M}_P(f,g)
-\frac{\mathsf{M}_P(f,\one)\mathsf{M}_P(g,\one)}
{\mathsf{M}_P(\one,\one)}=\mathsf{M}_P(f-\bar{f}\one,g-\bar{g}\one).
}

Assume now that the facet normals occur in antipodal pairs, indexed by a \emph{signed label set}
\eq{
I=I^{+}\sqcup(-I^{+}),
\quad
0\notin I,
\quad
u_{-i}=-u_i
\quad(i\in I^{+}).
}
The set $I^{+}$ contains one label from each pair. A support vector $h\in\R^I$ and a facet function $f\in\R^I$ are \emph{even} when
\eq{
h_{-i}=h_i,
\quad
f_{-i}=f_i
\quad(i\in I^{+}),
}
respectively. For a complete irredundant facet representation with signed label set $I$, the support vector $h$ is even if and only if $P=-P$.

\begin{definition}
The \emph{first nonconstant even eigenvalue} of a full-dimensional origin-symmetric polytope $P\subseteq\R^n$ is defined as
\eq{
\label{eq:sec-2-first-even-Rayleigh}
\lambda_{1,e}(P)=\min_{\substack{f\ \mathrm{even}\\ f\ \mathrm{not\ constant}}}
\frac{\mathsf{L}_P(f,f)}{\widetilde{\mathsf{M}}_P(f,f)}.
}
\end{definition}

For prescribed antipodal unit normals $(u_i)_{i\in I}$ that are pairwise distinct and positively span $\R^n$, the associated half-space intersection is
\eq{
P(q)=\bigcap_{i\in I}\{x\in\R^n:\ip{x}{u_i}\leq q_i\},
\quad q\in(0,\infty)^I.
}
The irredundant support domain $\mathcal{U}_{\mathrm{irr}}^I$ consists of the vectors $q\in(0,\infty)^I$ for which each inequality $\ip{x}{u_i}\leq q_i$ in the representation of $P(q)$ defines a facet of $P(q)$; see \cite[Sec. 8]{HI26}. If $q\in\mathcal{U}_{\mathrm{irr}}^I$ is even, then $P(q)$ is origin-symmetric. An even support vector $h\in\mathcal{U}_{\mathrm{irr}}^I$ is a \emph{local minimizer} of $\lambda_{1,e}$ under even support variations if $\lambda_{1,e}(P(q))\geq \lambda_{1,e}(P(h))$ for every even $q\in\mathcal{U}_{\mathrm{irr}}^I$ sufficiently close to $h$; the normal fan need not be fixed.

For nonempty compact sets $K,L\subseteq\R^n$, their \emph{Hausdorff distance} is
\eq{
d_H(K,L)=\max\left\{
\sup_{x\in K}\inf_{y\in L}|x-y|,
\quad
\sup_{y\in L}\inf_{x\in K}|x-y|
\right\}.
}
For nonempty compact convex sets, the Hausdorff distance has the support-function representation \cite[p. 66, Lem. 1.8.14]{Sch14}
\eq{
d_H(K,L)=\max_{u\in\Sn}|h_K(u)-h_L(u)|.
}
Blaschke's selection theorem states that every sequence of nonempty compact convex sets contained in a fixed compact set has a Hausdorff-convergent subsequence \cite[p. 63, Thm. 1.8.7]{Sch14}.

For a convex body $K$, its \emph{surface-area measure} is the Borel measure on $\Sn$ defined by
\eq{
S_K(\Omega)=\mathcal{H}^{n-1}\bigl(\{x\in\partial K:\nu_K(x)\in\Omega\}\bigr)
}
for every Borel set $\Omega\subseteq\Sn$. Here $\nu_K(x)$ is the outer unit normal, which is unique for $\mathcal{H}^{n-1}$-almost every $x\in\partial K$. If $P$ is a full-dimensional polytope with facets $F_i$ and outer unit normals $u_i$, then
\eq{
S_P=\sum_i\mathcal{H}^{n-1}(F_i)\delta_{u_i},
}
where $\delta_{u_i}$ is the Dirac measure at $u_i$. A sequence of finite Borel measures $\mu_{\nu}$ on $\Sn$ converges weakly to $\mu$, written $\mu_{\nu}\rightharpoonup\mu$, when
\eq{
\int_{\Sn}\varphi\,d\mu_{\nu}\to \int_{\Sn}\varphi\,d\mu
}
for every continuous function $\varphi$ on $\Sn$.

For $k\geq 2$ (including $k=\infty$), an origin-symmetric convex body $K$ belongs to $C^k_{+}$ if $h_K\in C^k(\Sn)$ and $r_K=\hess h_K+h_K\bar{g}$ is positive definite. Here $\hess$ denotes the Riemannian Hessian with respect to the round metric $\bar{g}$, and $d\omega$ is the corresponding Riemannian measure.

For a convex body $K$ containing the origin in its interior, its \emph{Minkowski functional} is
\eq{
\norm{x}_K=\inf\{t> 0:x\in tK\}
\quad(x\in\R^n).
}
The Minkowski functional is convex and positively homogeneous of degree one, and satisfies
\eq{
K=\{x\in\R^n:\norm{x}_K\leq 1\},
\quad
\norm{x}_K=h_{K^{\ast}}(x)
\quad(x\in\R^n).
}

The following continuity properties are standard; see \cite[p. 68, Thm. 1.8.20, p. 212, Thm. 4.2.1]{Sch14}. If nonempty compact convex sets $K_{\nu}$ converge to $K$ in the Hausdorff metric, then $\Vol(K_{\nu})\to \Vol(K)$. If $K_{\nu}$ and $K$ are convex bodies, then $S_{K_{\nu}}\rightharpoonup S_K$, and every compact set $E\subseteq\operatorname{int}K$ is contained in $\operatorname{int}K_{\nu}$ for all sufficiently large $\nu$. If in addition there are constants $0<r\leq R$ such that $rB_2^n\subseteq K_{\nu}\subseteq RB_2^n$ for every $\nu$, then $rB_2^n\subseteq K\subseteq RB_2^n$, the Minkowski functionals satisfy
\eq{
\norm{\mathord{\cdot}}_{K_{\nu}}\to \norm{\mathord{\cdot}}_K
\quad\text{locally uniformly on }\R^n,
}
and $K_{\nu}^{\ast}\to K^{\ast}$ in the Hausdorff metric.

\section{Continuity at polyhedral limits}

\subsection{Radial extensions and the dual energy}

All derivatives in this section are Euclidean. For an open set $O\subseteq\R^n$, let $C_c^{\infty}(O)$ denote the smooth functions with compact support in $O$. A locally integrable function $g$ defines a distribution by
\eq{
\langle g,\varphi\rangle=\int_Og\varphi\,dx,
\quad \varphi\in C_c^{\infty}(O).
}
For a distribution $\mathfrak{u}$, its derivative $D_k\mathfrak{u}$ is defined by
\eq{
\langle D_k\mathfrak{u},\varphi\rangle=-\langle \mathfrak{u},\partial_k\varphi\rangle.
}
The \emph{distributional gradient} $Dg=(D_kg)_{k=1}^n$ and \emph{distributional Hessian} $D^2g=(D^2_{kl}g)_{k,l=1}^n$, with $D^2_{kl}g=D_l(D_kg)$, satisfy
\eq{
\label{eq:sec-3-distributional-definitions}
\langle D_kg,\varphi\rangle
&=-\int_Og\partial_k\varphi\,dx,\\
\langle D^2_{kl}g,\varphi\rangle
&=-\langle D_kg,\partial_l\varphi\rangle
=\int_Og\partial_k\partial_l\varphi\,dx.
}
For a test field $X=(X_1,\ldots,X_n)\in C_c^{\infty}(O;\R^n)$, we write
\eq{
\label{eq:sec-3-gradient-field-pairing}
\langle Dg,X\rangle
=\sum_{k=1}^n\langle D_kg,X_k\rangle
=-\int_Og\sum_{k=1}^n\partial_kX_k\,dx
=-\int_Og\operatorname{div}X\,dx.
}
A nonnegative Borel measure $\mu$ on $\R^n$ is \emph{locally finite} if $\mu(K)<\infty$ for every compact set $K\subseteq\R^n$. A signed, vector-valued, or matrix-valued Borel measure is locally finite if the total variation of each component is finite on every compact set. For a Borel measure $\mu$ and a Borel set $E$, the restriction of $\mu$ to $E$ is defined by $(\mu\mathbin{\llcorner}E)(G)=\mu(G\cap E)$ for Borel sets $G$.

Write $B_r(x_0)=\{x\in\R^n:|x-x_0|<r\}$ for $r>0$. Let $P$ be a full-dimensional polytope containing the origin in its interior, with facet data $(F_i,u_i,h_i)_{i\in J}$. For a nonempty convex set $F$, define $\operatorname{pos}F=\{tx:t\geq 0,\ x\in F\}$, and set $\operatorname{pos}\varnothing=\{0\}$. Put $v_i=u_i/h_i$ and $C_i=\operatorname{pos}F_i$ for $i\in J$. Since $P=\bigcap_{i\in J}\{x\in\R^n:\ip{v_i}{x}\leq 1\}$, the Minkowski functional satisfies
\eq{
\norm{x}_P=\max_{i\in J}\ip{v_i}{x}
\quad(x\in\R^n).
}
For $x\neq 0$, we have $\norm{x}_P> 0$ and $x/\norm{x}_P\in\partial P$. Therefore
\eq{
\label{eq:sec-3-cone-inequalities}
C_i
=\{x:\ip{v_i}{x}=\norm{x}_P\}
&=\bigcap_{k\in J}\{x:\ip{v_k-v_i}{x}\leq 0\},\\
C_i\cap\{x:\ip{v_i}{x}=1\}&=F_i,\quad
C_i\cap C_j=\operatorname{pos}(F_i\cap F_j).
}
The $C_i$ are full-dimensional polyhedral cones with disjoint interiors and union $\R^n$. The facets of $C_i$ are the cones over the ridges of $P$ contained in $F_i$.

Let $\mathsf{E}$ be the set of unordered adjacent facet pairs, and put
\eq{
C_{ij}=\operatorname{pos}R_{ij},
\quad
a_{ij}=|v_i-v_j|.
}
For $\{i,j\}\in\mathsf{E}$, we have $a_{ij}> 0$ because $u_i\neq u_j$ and $h_i,h_j> 0$, and we call the common closed codimension-one face $C_{ij}=C_i\cap C_j$ the \emph{interface} between $C_i$ and $C_j$. The relative boundaries of the interfaces $C_{ij}$, as well as the intersections of distinct interfaces, are finite unions of cones of dimension at most $n-2$. In particular, for
\eq{
\mathcal{Z}=\bigcup_{\{i,j\}\in\mathsf{E}}
\operatorname{relbd}C_{ij},
\quad
\mathcal{S}=\bigcup_{\{i,j\}\in\mathsf{E}}C_{ij},
}
we have $\mathcal{H}^{n-1}(\mathcal{Z})=0$ and $\mathcal{H}^n(\mathcal{S})=0$.

For each ordered pair $(i,j)$ with $\{i,j\}\in\mathsf{E}$, define
\eq{
\eta_{ij}=\frac{v_j-v_i}{a_{ij}},
}
which satisfies $\eta_{ji}=-\eta_{ij}$. Since
\eq{
C_i\subseteq\{x:\ip{\eta_{ij}}{x}\leq 0\},
\quad
C_j\subseteq\{x:\ip{\eta_{ij}}{x}\geq 0\},
\quad
C_{ij}\subseteq\eta_{ij}^{\perp},
}
the vector $\eta_{ij}$ is the outward unit normal of $C_i$ on $\relint C_{ij}$, and $\eta_{ji}=-\eta_{ij}$ is the outward unit normal of $C_j$ on $\relint C_{ij}$. Every sum over $\{i,j\}\in\mathsf{E}$ counts each unordered adjacent pair once with either order $(i,j)$ or $(j,i)$.

For a facet function $f$, its \emph{radial extension} $\hat{f}$ is defined by
\eq{
\hat{f}(x)=f_i\quad(x\in\operatorname{int}C_i).
}
The values of $\hat{f}$ on $\mathcal{S}$ do not affect any Lebesgue integral. Note that $\hat{f}\in L^2_{\mathrm{loc}}(\R^n)$ and $\hat{f}(tx)=\hat{f}(x)$ for every $t> 0$ and almost every $x\in\R^n$.

\begin{lemma}
\label{lem:sec-3-distributional-derivatives}
The distributional gradient of $\norm{\mathord{\cdot}}_P$ is
\eq{
\label{eq:sec-3-gauge-gradient}
D\norm{\mathord{\cdot}}_P=\mathbf{v}_P\,dx,
\quad
\mathbf{v}_P=\sum_{i\in J}v_i\mathbf{1}_{\operatorname{int}C_i}.
}
For every facet function $f$ and its radial extension $\hat{f}$, the distributional gradient of $\hat{f}$ and Hessian of $\norm{\mathord{\cdot}}_P$ are the locally finite measures
\eq{
\label{eq:sec-3-distributional-derivatives}
D\hat{f}&=\sum_{\{i,j\}\in\mathsf{E}}
(f_j-f_i)\eta_{ij}
\mathcal{H}^{n-1}\mathbin{\llcorner}C_{ij},\\
D^2\norm{\mathord{\cdot}}_P&=\sum_{\{i,j\}\in\mathsf{E}}
a_{ij}\eta_{ij}\otimes\eta_{ij}
\mathcal{H}^{n-1}\mathbin{\llcorner}C_{ij}.
}
\end{lemma}

\begin{proof}
Let $\varphi\in C_c^{\infty}(\R^n)$, choose $r>0$ with $\operatorname{supp}\varphi\subseteq Q=(-r,r)^n$, and put $K_i=C_i\cap\overline{Q}$. Write $\nu_i^Q$ for the outward unit normal of $K_i$, defined almost everywhere on $\partial K_i$. For each coordinate vector $\mathbf{e}_l$, we apply the Gauss--Green formula on the polyhedral domain $\operatorname{int}K_i$ \cite[pp. 92--93, Rem. 9.5]{Mag12} to obtain
\eq{
\label{eq:sec-3-cone-integration-by-parts}
\int_{C_i}\partial_l\varphi\,dx
=\int_{K_i}\operatorname{div}(\varphi\mathbf{e}_l)\,dx
=\int_{\partial K_i}\varphi(\nu_i^Q)_l\,d\mathcal{H}^{n-1}=\sum_{j:\{i,j\}\in\mathsf{E}}
(\eta_{ij})_l\int_{C_{ij}}\varphi\,d\mathcal{H}^{n-1}.
}
Note that here we use that the facets of $K_i$ contained in $\partial Q$ contribute zero as $\varphi$ vanishes near $\partial Q$. Moreover, every remaining facet of $K_i$ is a truncated interface $C_{ij}\cap\overline{Q}$ and the faces of $K_i$ of dimension at most $n-2$ have zero $\mathcal{H}^{n-1}$-measure.

The first identity in \eqref{eq:sec-3-distributional-derivatives} now follows by applying \eqref{eq:sec-3-cone-integration-by-parts} to $\langle D_l\hat{f},\varphi\rangle$:
\eq{
\langle D_l\hat{f},\varphi\rangle
&=-\int_{\R^n}\hat{f}\partial_l\varphi\,dx
=-\sum_{i\in J}\int_{\operatorname{int}C_i}f_i\partial_l\varphi\,dx\\
&=-\sum_{i\in J}f_i\int_{C_i}\partial_l\varphi\,dx
=-\sum_{i\in J}\sum_{j:\{i,j\}\in\mathsf{E}}
f_i(\eta_{ij})_l\int_{C_{ij}}\varphi\,d\mathcal{H}^{n-1}\\
&=\sum_{\{i,j\}\in\mathsf{E}}
(-f_i(\eta_{ij})_l-f_j(\eta_{ji})_l)
\int_{C_{ij}}\varphi\,d\mathcal{H}^{n-1}\\
&=\sum_{\{i,j\}\in\mathsf{E}}
(f_j-f_i)(\eta_{ij})_l
\int_{C_{ij}}\varphi\,d\mathcal{H}^{n-1}.
}

Next we compute the gradient of $\norm{\mathord{\cdot}}_P$. For each $i\in J$, we have
\eq{
\partial_k\bigl(\ip{v_i}{x}\varphi\bigr)
=(v_i)_k\varphi+\ip{v_i}{x}\partial_k\varphi.
}
Applying \eqref{eq:sec-3-cone-integration-by-parts} to $\ip{v_i}{x}\varphi$ and rearranging, we obtain
\eq{
\int_{C_i}\ip{v_i}{x}\partial_k\varphi\,dx
=\sum_{j:\{i,j\}\in\mathsf{E}}(\eta_{ij})_k
\int_{C_{ij}}\ip{v_i}{x}\varphi\,d\mathcal{H}^{n-1}
-(v_i)_k\int_{C_i}\varphi\,dx.
}
Since $\norm{x}_P=\ip{v_i}{x}$ on $C_i$ and $\ip{v_i-v_j}{x}=0$ on $C_{ij}$, summing over $i\in J$ yields
\eq{
\label{eq:sec-3-gauge-gradient-pairing}
\langle D_k\norm{\mathord{\cdot}}_P,\varphi\rangle
&=-\sum_{i\in J}\int_{C_i}\ip{v_i}{x}\partial_k\varphi\,dx\\
&=\sum_{i\in J}(v_i)_k\int_{C_i}\varphi\,dx
-\sum_{i\in J}\sum_{j:\{i,j\}\in\mathsf{E}}(\eta_{ij})_k
\int_{C_{ij}}\ip{v_i}{x}\varphi\,d\mathcal{H}^{n-1}\\
&=\sum_{i\in J}(v_i)_k\int_{C_i}\varphi\,dx
-\sum_{\{i,j\}\in\mathsf{E}}(\eta_{ij})_k
\int_{C_{ij}}\ip{v_i-v_j}{x}\varphi\,d\mathcal{H}^{n-1}\\
&=\sum_{i\in J}(v_i)_k\int_{C_i}\varphi\,dx.
}

By \eqref{eq:sec-3-gauge-gradient-pairing} and \eqref{eq:sec-3-cone-integration-by-parts}, we obtain
\eq{
\langle D^2_{kl}\norm{\mathord{\cdot}}_P,\varphi\rangle
&=-\langle D_k\norm{\mathord{\cdot}}_P,\partial_l\varphi\rangle
=-\sum_{i\in J}(v_i)_k\int_{C_i}\partial_l\varphi\,dx\\
&=\sum_{\{i,j\}\in\mathsf{E}}
(v_j-v_i)_k(\eta_{ij})_l
\int_{C_{ij}}\varphi\,d\mathcal{H}^{n-1}\\
&=\sum_{\{i,j\}\in\mathsf{E}}
a_{ij}(\eta_{ij})_k(\eta_{ij})_l
\int_{C_{ij}}\varphi\,d\mathcal{H}^{n-1}.
}
Since $\varphi$ is arbitrary, the second identity in \eqref{eq:sec-3-distributional-derivatives} follows.
\end{proof}

\begin{definition}
For a full-dimensional polytope $P\subseteq\R^n$ containing the origin in its interior and a bounded open set $\Omega\subseteq\R^n$, we define the \emph{dual energy} of $g\in L^2_{\mathrm{loc}}(\R^n)$ by
\eq{
\label{eq:sec-3-dual-energy}
\mathcal{L}_{\Omega}(P,g)
=\sup_{X\in C_c^{\infty}(\Omega;\R^n)}
\left(-2\int_{\Omega}g\operatorname{div}X\,dx
-\int_{\Omega}X^{\mathsf{T}}\,d(D^2\norm{\mathord{\cdot}}_P)X\right),
}
where
\eq{
\int_{\Omega}X^{\mathsf{T}}\,d(D^2\norm{\mathord{\cdot}}_P)X
=\sum_{k,l=1}^n\int_{\Omega}X_kX_l\,d(D^2_{kl}\norm{\mathord{\cdot}}_P).
}
\end{definition}

By \autoref{lem:sec-3-distributional-derivatives}, each entry $D^2_{kl}\norm{\mathord{\cdot}}_P$ is a locally finite signed Borel measure.

\begin{theorem}
For every facet function $f$ on $P$ and every bounded open set $\Omega\subseteq\R^n$, the radial extension $\hat{f}$ satisfies
\eq{
\label{eq:sec-3-interface-energy}
\mathcal{L}_{\Omega}(P,\hat{f})
=\sum_{\{i,j\}\in\mathsf{E}}
\frac{(f_i-f_j)^2}{a_{ij}}
\mathcal{H}^{n-1}(C_{ij}\cap\Omega).
}
Moreover,
\eq{
\label{eq:sec-3-radial-mass}
\mathsf{M}_P(f,f)=n\int_P\hat{f}^2\,dx,
\quad
\mathsf{M}_P(f,\one)=n\int_P\hat{f}\,dx,
}
and
\eq{
\label{eq:sec-3-radial-energy}
\mathsf{L}_P(f,f)=(n-1)\mathcal{L}_{\operatorname{int}P}(P,\hat{f}).
}
\end{theorem}

\begin{proof}
Let us write $b_{ij}=f_j-f_i$ and fix $X\in C_c^{\infty}(\Omega;\R^n)$, extended by zero outside $\Omega$. It follows from \eqref{eq:sec-3-gradient-field-pairing} and \autoref{lem:sec-3-distributional-derivatives} that
\eq{
-\int_{\Omega}\hat{f}\operatorname{div}X\,dx
&=\sum_{k=1}^n\langle D_k\hat{f},X_k\rangle\\
&=\sum_{k=1}^n\sum_{\{i,j\}\in\mathsf{E}}
b_{ij}(\eta_{ij})_k
\int_{C_{ij}\cap\Omega}X_k\,d\mathcal{H}^{n-1}\\
&=\sum_{\{i,j\}\in\mathsf{E}}b_{ij}
\int_{C_{ij}\cap\Omega}\ip{X}{\eta_{ij}}\,d\mathcal{H}^{n-1}.
}
Similarly,
\eq{
\int_{\Omega}X^{\mathsf{T}}\,d(D^2\norm{\mathord{\cdot}}_P)X
&=\sum_{k,l=1}^n\sum_{\{i,j\}\in\mathsf{E}}
a_{ij}(\eta_{ij})_k(\eta_{ij})_l
\int_{C_{ij}\cap\Omega}X_kX_l\,d\mathcal{H}^{n-1}\\
&=\sum_{\{i,j\}\in\mathsf{E}}a_{ij}
\int_{C_{ij}\cap\Omega}
\left(\sum_{k=1}^n(\eta_{ij})_kX_k\right)^2\,d\mathcal{H}^{n-1}\\
&=\sum_{\{i,j\}\in\mathsf{E}}a_{ij}
\int_{C_{ij}\cap\Omega}\ip{X}{\eta_{ij}}^2\,d\mathcal{H}^{n-1}.
}
Therefore, the expression under the supremum in \eqref{eq:sec-3-dual-energy} equals
\eq{
\sum_{\{i,j\}\in\mathsf{E}}
\int_{C_{ij}\cap\Omega}
\bigl(2b_{ij}\ip{X}{\eta_{ij}}-a_{ij}\ip{X}{\eta_{ij}}^2\bigr)
\,d\mathcal{H}^{n-1}.
}
Note that
\eq{
2b_{ij}\ip{X}{\eta_{ij}}-a_{ij}\ip{X}{\eta_{ij}}^2
=\frac{b_{ij}^2}{a_{ij}}
-a_{ij}\left(\ip{X}{\eta_{ij}}-\frac{b_{ij}}{a_{ij}}\right)^2\leq \frac{b_{ij}^2}{a_{ij}}.
}
Therefore,
\eq{
\mathcal{L}_{\Omega}(P,\hat{f})
\leq \sum_{\{i,j\}\in\mathsf{E}}
\frac{b_{ij}^2}{a_{ij}}
\mathcal{H}^{n-1}(C_{ij}\cap\Omega).
}

For the reverse direction, fix $\varepsilon> 0$. By inner regularity and $\mathcal{H}^{n-1}(\mathcal{Z})=0$, we can choose compact sets $K_{ij}\subseteq\Omega\cap\relint C_{ij}$ such that
\eq{
\sum_{\{i,j\}\in\mathsf{E}}\frac{b_{ij}^2}{a_{ij}}
\mathcal{H}^{n-1}\bigl((C_{ij}\cap\Omega)\setminus K_{ij}\bigr)
<\varepsilon.
}
Note that when $K_{ij}\neq \varnothing$, $K_{ij}$ has positive distance from $\Omega^c$ and from each $C_{kl}$ with $\{k,l\}\neq \{i,j\}$. In particular, we can choose pairwise disjoint open neighborhoods $V_{ij}$ of the sets $K_{ij}$, with $V_{ij}\Subset\Omega$ and $V_{ij}\cap C_{kl}=\varnothing$ for $\{k,l\}\neq \{i,j\}$. 

For $K_{ij}=\varnothing$, set $\chi_{ij}=0$; for $K_{ij}\neq \varnothing$, take $\chi_{ij}\in C_c^{\infty}(V_{ij})$ with $0\leq \chi_{ij}\leq 1$ and $\chi_{ij}=1$ on $K_{ij}$. Define $X=\sum_{\{i,j\}\in\mathsf{E}} \frac{b_{ij}}{a_{ij}}\chi_{ij}\eta_{ij}$. Then $X\in C_c^{\infty}(\Omega;\R^n)$, and all cutoffs except $\chi_{ij}$ vanish on $C_{ij}$. Since $2\chi_{ij}-\chi_{ij}^2\geq 0$ everywhere and equals one on $K_{ij}$, it follows that
\eq{
\mathcal{L}_{\Omega}(P,\hat{f})
&\geq \sum_{\{i,j\}\in\mathsf{E}}\frac{b_{ij}^2}{a_{ij}}
\int_{C_{ij}\cap\Omega}
(2\chi_{ij}-\chi_{ij}^2)\,d\mathcal{H}^{n-1}\\
&\geq \sum_{\{i,j\}\in\mathsf{E}}\frac{b_{ij}^2}{a_{ij}}
\mathcal{H}^{n-1}(K_{ij})\\
&\geq \sum_{\{i,j\}\in\mathsf{E}}\frac{b_{ij}^2}{a_{ij}}
\mathcal{H}^{n-1}(C_{ij}\cap\Omega)-\varepsilon.
}
Letting $\varepsilon\to 0^{+}$ establishes the reverse inequality.

Next we prove \eqref{eq:sec-3-radial-mass} and \eqref{eq:sec-3-radial-energy}. Let $V\subseteq\R^n$ be an $m$-dimensional linear subspace with $m\geq 2$, and let $B\subseteq V$ be a compact convex set in $\{x\in V:\ip{x}{u}=d\}$ for a unit vector $u\in V$ and $d> 0$. For each $s\in(0,d)$, the slice $\{x\in\operatorname{conv}(\{0\}\cup B):\ip{x}{u}=s\}$ equals $(s/d)B$. By Fubini's theorem,
\eq{
\label{eq:sec-3-pyramid-volume}
\mathcal{H}^m\bigl(\operatorname{conv}(\{0\}\cup B)\bigr)
=\int_0^d\left(\frac{s}{d}\right)^{m-1}\mathcal{H}^{m-1}(B)\,ds
=\frac{d}{m}\mathcal{H}^{m-1}(B).
}

By \eqref{eq:sec-3-cone-inequalities}, for each $i\in J$, we have $C_i\cap P=\operatorname{conv}(\{0\}\cup F_i)$ with $F_i\subseteq\{x\in\R^n:\ip{x}{u_i}=h_i\}$. Now applying \eqref{eq:sec-3-pyramid-volume} yields
\eq{
\Vol(C_i\cap P)
=\frac{h_i}{n}\mathcal{H}^{n-1}(F_i)=\frac{m_i}{n}.
}
The pyramids $C_i\cap P$ cover $P$ with pairwise disjoint interiors, and their boundaries have zero volume. Since $\hat{f}=f_i$ on $\operatorname{int}C_i$, it follows that
\eq{
n\int_P\hat{f}^2\,dx
&=n\sum_{i\in J}f_i^2\Vol(C_i\cap P)
=\sum_{i\in J}m_if_i^2,\\
n\int_P\hat{f}\,dx
&=n\sum_{i\in J}f_i\Vol(C_i\cap P)
=\sum_{i\in J}m_if_i.
}

To obtain \eqref{eq:sec-3-radial-energy}, we need to compute the measure of each interface inside $P$, namely $C_{ij}\cap\operatorname{int}P$, and relate it to the ridge conductance. For an adjacent pair $\{i,j\}\in\mathsf{E}$, put
\eq{
H_{ij}=\eta_{ij}^{\perp},\quad
w_{ij}=v_i-\ip{v_i}{\eta_{ij}}\eta_{ij}\in H_{ij}.
}
The affine hull of $R_{ij}$ is the intersection of the supporting hyperplanes of $F_i$ and $F_j$:
\eq{
\operatorname{aff}R_{ij}
=\{x:\ip{v_i}{x}=\ip{v_j}{x}=1\}=\{x\in H_{ij}:\ip{w_{ij}}{x}=1\}.
}
Therefore, $d_{ij}:=\operatorname{dist}(0,\operatorname{aff}R_{ij})=\frac{1}{|w_{ij}|}$. Moreover,
\eq{
\label{eq:sec-3-ridge-distance}
d_{ij}=\frac{a_{ij}}{|v_i\wedge v_j|},
\quad
|v_i\wedge v_j|=\frac{\sin\theta_{ij}}{h_ih_j}.
}

On the other hand, by $\norm{x}_P=\ip{v_i}{x}$ on $C_{ij}$, we have
\eq{
C_{ij}\cap P=\operatorname{conv}(\{0\}\cup R_{ij}),
\quad
(C_{ij}\cap P)\setminus(C_{ij}\cap\operatorname{int}P)=R_{ij},\quad
\mathcal{H}^{n-1}(R_{ij})=0.
}
Now applying \eqref{eq:sec-3-pyramid-volume} in $H_{ij}$ with $m=n-1$ yields
\eq{
\mathcal{H}^{n-1}(C_{ij}\cap\operatorname{int}P)
=\mathcal{H}^{n-1}(C_{ij}\cap P)
=\frac{d_{ij}}{n-1}\mathcal{H}^{n-2}(R_{ij}).
}
Therefore, by \eqref{eq:sec-3-ridge-distance},
\eq{
\frac{n-1}{a_{ij}}
\mathcal{H}^{n-1}(C_{ij}\cap\operatorname{int}P)
=\frac{d_{ij}}{a_{ij}}\mathcal{H}^{n-2}(R_{ij})
=\frac{h_ih_j\mathcal{H}^{n-2}(R_{ij})}{\sin\theta_{ij}}
=c_{ij}.
}
Substituting into \eqref{eq:sec-3-interface-energy}, we obtain
\eq{
(n-1)\mathcal{L}_{\operatorname{int}P}(P,\hat{f})
=\sum_{\{i,j\}\in\mathsf{E}}c_{ij}(f_j-f_i)^2
=\mathsf{L}_P(f,f).
}
\end{proof}

\subsection{\texorpdfstring{Local $L^2$ compactness of radial extensions}{Local L2 compactness of radial extensions}}

For an open set $O\subseteq\R^n$ and $p\in\{1,2\}$, the space $L^p(O)$ consists of measurable functions, identified up to equality almost everywhere, with finite norm
\eq{
\norm{w}_{L^p(O)}=\left(\int_O|w(x)|^p\,dx\right)^{\frac{1}{p}}.
}
A sequence $w_{\nu}$ converges to $w$ in the $L^p(O)$-norm when $\norm{w_{\nu}-w}_{L^p(O)}\to 0$. A measurable function $w$ belongs to $L^p_{\mathrm{loc}}(\R^n)$ when $w|_{B_r(0)}\in L^p(B_r(0))$ for every $r>0$. We write $w_{\nu}\to w$ in $L^p_{\mathrm{loc}}(\R^n)$ when $w_{\nu}\to w$ in the $L^p(B_r(0))$-norm for every $r>0$.

A function $w\in L^1(O)$ is of bounded variation, written $w\in\mathrm{BV}(O)$, if the distributional gradient $Dw$ is a finite vector-valued Borel measure. The \emph{total variation measure} $|Dw|$ is characterized, for every open set $U\subseteq O$, by
\eq{
|Dw|(U)
=\sup_{\substack{X\in C_c^{\infty}(U;\R^n)\\|X(x)|\leq 1}}
\left(-\int_U w\operatorname{div}X\,dx\right).
}
Set $\norm{w}_{\mathrm{BV}(O)}=\norm{w}_{L^1(O)}+|Dw|(O)$. By the compactness theorem for functions of bounded variation \cite[p. 17, Thm. 1.19]{Giu84}, every sequence bounded in $\mathrm{BV}(O)$ on a bounded open set $O\subseteq\R^n$ with Lipschitz boundary has a subsequence converging in $L^1(O)$.

For a facet function $f=(f_i)_{i\in J}$, the product $f\lvert f\rvert$ is understood componentwise, namely, $f\lvert f\rvert=(f_i\lvert f_i\rvert)_{i\in J}$; its radial extension is $\widehat{f\lvert f\rvert}=\hat{f}\lvert\hat{f}\rvert$.

\begin{lemma}\label{lem:sec-3-radial-compactness}
Suppose that $P\subseteq\R^n$ is a polytope with $B_2^n\subseteq P$, let $f$ be a facet function on $P$, and put $w=\hat{f}\lvert\hat{f}\rvert$. Then for every $r>0$, we have $w\in\mathrm{BV}(B_r(0))$ and
\eq{
\label{eq:sec-3-signed-square-BV-bound}
\norm{w}_{L^1(B_r(0))}
&\leq \frac{r^n}{n}\mathsf{M}_P(f,f),\\
|Dw|(B_r(0))
&\leq r^{n-1}\left(\mathsf{M}_P(f,f)
+\frac{4}{n-1}\mathsf{L}_P(f,f)\right).
}
Consequently, if polytopes $P_{\nu}\subseteq\R^n$ satisfy $B_2^n\subseteq P_{\nu}$ and facet functions $f_{\nu}$ on $P_{\nu}$ satisfy
\eq{
\mathsf{M}_{P_{\nu}}(f_{\nu},f_{\nu})
+\mathsf{L}_{P_{\nu}}(f_{\nu},f_{\nu})\leq C_0
}
for some constant $C_0$, then there exist a subsequence and a function $g\in L^2_{\mathrm{loc}}(\R^n)$ such that
\eq{
\norm{\hat{f}_{\nu}-g}_{L^2(B_r(0))}\to 0
\quad\text{for every }r>0.
}
\end{lemma}

\begin{proof}
Since $B_2^n\subseteq P$, we have $h_i\geq 1$ and $|v_i|=1/h_i\leq 1$. Set $W=\hat{f}^2\mathbf{v}_P$, where $\mathbf{v}_P=v_i$ on $\operatorname{int}C_i$ as in \eqref{eq:sec-3-gauge-gradient}. On $\operatorname{int}C_i$, we have $w=f_i\lvert f_i\rvert$ and $W=f_i^2v_i$.

For each ordered pair $(i,j)$ with $\{i,j\}\in\mathsf{E}$, set
\eq{
\mu_{ij}=\mathcal{H}^{n-1}\mathbin{\llcorner}C_{ij},\quad
\beta_{ij}=\ip{f_j^2v_j-f_i^2v_i}{\eta_{ij}},
}
and define the nonnegative, locally finite measure
\eq{
\mathcal{E}=\sum_{\{i,j\}\in\mathsf{E}}
\frac{(f_j-f_i)^2}{a_{ij}}\mu_{ij}.
}
By \eqref{eq:sec-3-interface-energy}, $\mathcal{E}(\Omega)=\mathcal{L}_{\Omega}(P,\hat{f})$ for every bounded open set $\Omega$. We will show that
\eq{
\label{eq:sec-3-measure-domination}
|Dw|\leq \operatorname{div}W+4\mathcal{E},
}
and then evaluate \eqref{eq:sec-3-measure-domination} on $B_r(0)$.

By \autoref{lem:sec-3-distributional-derivatives},
\eq{
Dw=\sum_{\{i,j\}\in\mathsf{E}}
\bigl(f_j\lvert f_j\rvert-f_i\lvert f_i\rvert\bigr)
\eta_{ij}\mu_{ij}.
}
The measures $\mu_{ij}$ are locally finite and have pairwise disjoint supports up to sets of zero $\mathcal{H}^{n-1}$-measure. Since $|\eta_{ij}|=1$, taking total variation yields
\eq{
\label{eq:sec-3-signed-square-variation}
|Dw|=\sum_{\{i,j\}\in\mathsf{E}}
\left|f_j\lvert f_j\rvert-f_i\lvert f_i\rvert\right|\mu_{ij}.
}
Note that for every $r>0$, $w$ takes only finitely many values almost everywhere and thus $w\in L^1(B_r(0))$. Moreover, by \eqref{eq:sec-3-signed-square-variation} and local finiteness of the measures $\mu_{ij}$, we have $|Dw|(B_r(0))<\infty$. Hence $w\in\mathrm{BV}(B_r(0))$.

Next we calculate $\operatorname{div}W$. Applying \autoref{lem:sec-3-distributional-derivatives} to the radial extension $W_k$ of $(f_i^2(v_i)_k)_{i\in J}$ for $k=1,\ldots,n$ yields
\eq{
D_kW_k
=\sum_{\{i,j\}\in\mathsf{E}}
\bigl(f_j^2(v_j)_k-f_i^2(v_i)_k\bigr)
(\eta_{ij})_k\mu_{ij}.
}
Therefore,
\eq{
\label{eq:sec-3-square-flux-divergence}
\operatorname{div}W
&=\sum_{k=1}^nD_kW_k
=\sum_{k=1}^n\sum_{\{i,j\}\in\mathsf{E}}
\bigl(f_j^2(v_j)_k-f_i^2(v_i)_k\bigr)
(\eta_{ij})_k\mu_{ij}\\
&=\sum_{\{i,j\}\in\mathsf{E}}
\left(\sum_{k=1}^n
\bigl(f_j^2(v_j)_k-f_i^2(v_i)_k\bigr)(\eta_{ij})_k\right)
\mu_{ij}\\
&=\sum_{\{i,j\}\in\mathsf{E}}
\ip{f_j^2v_j-f_i^2v_i}{\eta_{ij}}\mu_{ij}
=\sum_{\{i,j\}\in\mathsf{E}}\beta_{ij}\mu_{ij}.
}

We now compare the interface coefficients of $|Dw|$ and $\operatorname{div}W$. For each $\{i,j\}\in\mathsf{E}$, since $|v_i|,|v_j|\leq 1$ and $v_j-v_i=a_{ij}\eta_{ij}$,
\eq{
\beta_{ij}
&=\frac{a_{ij}}{2}(f_i^2+f_j^2)
+\frac{1}{2}(f_j-f_i)(f_i+f_j)\ip{v_i+v_j}{\eta_{ij}}\\
&\geq \frac{a_{ij}}{4}(|f_i|+|f_j|)^2
-|f_j-f_i|(|f_i|+|f_j|).
}
On the other hand,
\eq{
\left|f_j\lvert f_j\rvert-f_i\lvert f_i\rvert\right|
\leq |f_j-f_i|(|f_i|+|f_j|).
}
Combining the bounds for $\beta_{ij}$ and $\left|f_j\lvert f_j\rvert-f_i\lvert f_i\rvert\right|$, we obtain
\eq{
\label{eq:sec-3-signed-square-interface-bound}
&\beta_{ij}+\frac{4(f_j-f_i)^2}{a_{ij}}
-\left|f_j\lvert f_j\rvert-f_i\lvert f_i\rvert\right|\\
&\quad\geq \frac{a_{ij}}{4}(|f_i|+|f_j|)^2
-2|f_j-f_i|(|f_i|+|f_j|)
+\frac{4(f_j-f_i)^2}{a_{ij}}\\
&\quad=\left(\frac{\sqrt{a_{ij}}}{2}(|f_i|+|f_j|)
-\frac{2|f_j-f_i|}{\sqrt{a_{ij}}}\right)^2\geq 0.
}

In view of \eqref{eq:sec-3-signed-square-variation}, \eqref{eq:sec-3-square-flux-divergence}, and \eqref{eq:sec-3-signed-square-interface-bound}, we have
\eq{
\operatorname{div}W+4\mathcal{E}-|Dw|=\sum_{\{i,j\}\in\mathsf{E}}
\left(\beta_{ij}+\frac{4(f_j-f_i)^2}{a_{ij}}
-\left|f_j\lvert f_j\rvert-f_i\lvert f_i\rvert\right|\right)\mu_{ij}\geq 0.
}

Evaluating \eqref{eq:sec-3-measure-domination} on $B_r(0)$ and using $\mathcal{E}(B_r(0))=\mathcal{L}_{B_r(0)}(P,\hat{f})$, we obtain
\eq{
\label{eq:sec-3-signed-square-ball-bound}
|Dw|(B_r(0))
\leq (\operatorname{div}W)(B_r(0))+4\mathcal{L}_{B_r(0)}(P,\hat{f}).
}

The next step is to estimate $(\operatorname{div}W)(B_r(0))$. Up to sets of zero $\mathcal{H}^{n-1}$-measure,
\eq{
\partial(\operatorname{int}C_i\cap B_r(0))
=(\operatorname{int}C_i\cap\partial B_r(0))\cup\bigcup_{j:\{i,j\}\in\mathsf{E}}
(\relint C_{ij}\cap B_r(0)).
}
The omitted sets are contained in the union of $\mathcal{Z}$ and the intersections $C_{ij}\cap\partial B_r(0)$. Each intersection lies in the $(n-2)$-dimensional sphere $\eta_{ij}^{\perp}\cap\partial B_r(0)$ and has zero $\mathcal{H}^{n-1}$-measure. Hence $\operatorname{int}C_i\cap B_r(0)$ has almost $C^1$ boundary, with outward normals $x/r$ on the spherical piece $\operatorname{int}C_i\cap\partial B_r(0)$ and $\eta_{ij}$ on each planar piece $\relint C_{ij}\cap B_r(0)$. Choose $\chi\in C_c^{\infty}(\R^n)$ with $\chi=1$ on a neighborhood of $\overline{B_r(0)}$. Then $\chi=1$ and $\nabla\chi=0$ on a neighborhood of $\overline{\operatorname{int}C_i\cap B_r(0)}$. Applying the Gauss--Green theorem \cite[p. 93, Thm. 9.6]{Mag12} to $\chi f_i^2v_i$ on $\operatorname{int}C_i\cap B_r(0)$ and using $\mathcal{H}^{n-1}(\mathcal{Z})=0$, we obtain
\eq{
0&=\int_{\operatorname{int}C_i\cap B_r(0)}
\sum_{k=1}^n\partial_k\bigl(\chi f_i^2(v_i)_k\bigr)\,dx\\
&=\int_{\partial B_r(0)\cap\operatorname{int}C_i}
\chi f_i^2\ip{v_i}{\frac{x}{r}}\,d\mathcal{H}^{n-1}(x)\\
&\quad+\sum_{j:\{i,j\}\in\mathsf{E}}
\int_{C_{ij}\cap B_r(0)}
\chi f_i^2\ip{v_i}{\eta_{ij}}\,d\mathcal{H}^{n-1}(x)\\
&=\int_{\partial B_r(0)\cap\operatorname{int}C_i}
f_i^2\ip{v_i}{\frac{x}{r}}\,d\mathcal{H}^{n-1}(x)\\
&\quad+\sum_{j:\{i,j\}\in\mathsf{E}}
\int_{C_{ij}\cap B_r(0)}
f_i^2\ip{v_i}{\eta_{ij}}\,d\mathcal{H}^{n-1}(x)\\
&=\int_{\partial B_r(0)\cap\operatorname{int}C_i}
f_i^2\ip{v_i}{\frac{x}{r}}\,d\mathcal{H}^{n-1}(x)
+\sum_{j:\{i,j\}\in\mathsf{E}}
f_i^2\ip{v_i}{\eta_{ij}}\mu_{ij}(B_r(0)).
}
For each $\{i,j\}\in\mathsf{E}$, since $\eta_{ji}=-\eta_{ij}$, we have
\eq{
f_i^2\ip{v_i}{\eta_{ij}}
+f_j^2\ip{v_j}{\eta_{ji}}
=\ip{f_i^2v_i-f_j^2v_j}{\eta_{ij}}
=-\beta_{ij}.
}
Summing over $i\in J$, we obtain from \eqref{eq:sec-3-square-flux-divergence} that
\eq{
\label{eq:sec-3-divergence-boundary-bound}
(\operatorname{div}W)(B_r(0))
&=\sum_{\{i,j\}\in\mathsf{E}}\beta_{ij}\mu_{ij}(B_r(0))\\
&=\int_{\partial B_r(0)}\hat{f}^2\ip{\mathbf{v}_P}{\frac{x}{r}}\,d\mathcal{H}^{n-1}(x)\\
&\leq \int_{\partial B_r(0)}\hat{f}^2\,d\mathcal{H}^{n-1},
}
where we used $\ip{\mathbf{v}_P}{x/r}\leq 1$ on $\partial B_r(0)$. Moreover, we have
\eq{
\label{eq:sec-3-polar-mass-identity}
\int_{B_r(0)}\hat{f}^2\,dx
=\frac{r^n}{n}\int_{\Sn}\hat{f}(\theta)^2\,d\mathcal{H}^{n-1}(\theta)
=\frac{r}{n}\int_{\partial B_r(0)}\hat{f}^2\,d\mathcal{H}^{n-1}.
}
Combining \eqref{eq:sec-3-divergence-boundary-bound}, \eqref{eq:sec-3-polar-mass-identity}, and \eqref{eq:sec-3-signed-square-ball-bound}, we obtain
\eq{
\label{eq:sec-3-variation-mass-energy-bound}
|Dw|(B_r(0))
\leq \frac{n}{r}\int_{B_r(0)}\hat{f}^2\,dx+4\mathcal{L}_{B_r(0)}(P,\hat{f}).
}

Since $B_1(0)\subseteq\operatorname{int}P$, zero-homogeneity and \eqref{eq:sec-3-radial-mass} imply
\eq{
\label{eq:sec-3-local-mass-bound}
\norm{w}_{L^1(B_r(0))}
=\int_{B_r(0)}\hat{f}^2\,dx
=r^n\int_{B_1(0)}\hat{f}^2\,dx
\leq \frac{r^n}{n}\mathsf{M}_P(f,f).
}
Since $C_{ij}$ is a cone, dilation by $r$ implies
\eq{
\mu_{ij}(B_r(0))
=\mathcal{H}^{n-1}\bigl(r(C_{ij}\cap B_1(0))\bigr)
=r^{n-1}\mu_{ij}(B_1(0)).
}
It follows from \eqref{eq:sec-3-interface-energy} and \eqref{eq:sec-3-radial-energy} that
\eq{
\label{eq:sec-3-local-energy-bound}
\mathcal{L}_{B_r(0)}(P,\hat{f})
=r^{n-1}\mathcal{L}_{B_1(0)}(P,\hat{f})
\leq \frac{r^{n-1}}{n-1}\mathsf{L}_P(f,f).
}
Substituting \eqref{eq:sec-3-local-mass-bound} and \eqref{eq:sec-3-local-energy-bound} into \eqref{eq:sec-3-variation-mass-energy-bound}, we obtain \eqref{eq:sec-3-signed-square-BV-bound}.

Finally, we prove compactness. Put $w_{\nu}=\hat{f}_{\nu}\lvert\hat{f}_{\nu}\rvert$. By \eqref{eq:sec-3-signed-square-BV-bound} and our assumption that $\mathsf{M}_{P_{\nu}}(f_{\nu},f_{\nu})+\mathsf{L}_{P_{\nu}}(f_{\nu},f_{\nu})\leq C_0$, we have
\eq{
\sup_{\nu}\norm{w_{\nu}}_{\mathrm{BV}(B_m(0))}<\infty
\quad(m\in\mathbb{N}).
}
In view of \cite[p. 17, Thm. 1.19]{Giu84}, we can apply the compactness theorem for functions of bounded variation on each $B_m(0)$ and obtain a subsequence and a function $z\in L^1_{\mathrm{loc}}(\R^n)$ such that $w_{\nu}\to z$ in the $L^1(B_r(0))$-norm for every $r>0$.

Define
\eq{
g(x)=\operatorname{sgn}(z(x))\sqrt{|z(x)|},
\quad\operatorname{sgn}(0)=0.
}
Since $g\lvert g\rvert=z$, we have $g\in L^2_{\mathrm{loc}}(\R^n)$, and since $|s-t|^2\leq 2\left|s\lvert s\rvert-t\lvert t\rvert\right|$ for all $s,t\in\R$, we obtain, for every $r>0$,
\eq{
\norm{\hat{f}_{\nu}-g}_{L^2(B_r(0))}^2
\leq 2\norm{w_{\nu}-z}_{L^1(B_r(0))}\to 0.
}
\end{proof}

\subsection{Polyhedral limits and lower semicontinuity}

The next proposition proves that the limit $g$ in \autoref{lem:sec-3-radial-compactness} is the radial extension of a facet function $f$ on $P$, and proves lower semicontinuity of the energy form.

\begin{proposition}\label{prop:sec-3-facet-function-compactness}
Suppose that polytopes $P_{\nu}\subseteq\R^n$ satisfy, for a fixed $R\geq 1$,
\eq{
B_2^n\subseteq P_{\nu}\subseteq RB_2^n,
\quad
P_{\nu}\to P\quad\text{in the Hausdorff metric},
}
where $P$ is a polytope. Let $f_{\nu}$ be facet functions on $P_{\nu}$ such that
\eq{
\mathsf{M}_{P_{\nu}}(f_{\nu},f_{\nu})+\mathsf{L}_{P_{\nu}}(f_{\nu},f_{\nu})\leq C_0
}
for a constant $C_0$. After passing to a subsequence, there is a facet function $f$ on $P$ such that
\eq{
\label{eq:sec-3-radial-L2-convergence}
\norm{\hat{f}_{\nu}-\hat{f}}_{L^2(B_r(0))}\to 0
\quad\text{for every }r>0,\\
\mathsf{M}_{P_{\nu}}(f_{\nu},f_{\nu})\to \mathsf{M}_P(f,f),\quad
\mathsf{M}_{P_{\nu}}(f_{\nu},\one)\to \mathsf{M}_P(f,\one),
}
and
\eq{
\label{eq:sec-3-energy-liminf-o}
\mathsf{L}_P(f,f)\leq \liminf_{\nu\to \infty}\mathsf{L}_{P_{\nu}}(f_{\nu},f_{\nu}),
}
where $\liminf$ in \eqref{eq:sec-3-energy-liminf-o} is taken along the original sequence. Moreover, if every $f_{\nu}$ is centered, then $f$ is centered, and if $P_{\nu}=-P_{\nu}$ and every $f_{\nu}$ is even, then $f$ is even.
\end{proposition}

\begin{proof}
Choose a subsequence along which $\mathsf{L}_{P_{\nu}}(f_{\nu},f_{\nu})$ converges to the $\liminf$. By \autoref{lem:sec-3-radial-compactness}, we may pass to a further subsequence and find $g\in L^2_{\mathrm{loc}}(\R^n)$ such that
\eq{
\norm{\hat{f}_{\nu}-g}_{L^2(B_r(0))}\to 0
\quad\text{for every }r>0.
}
Since $B_2^n\subseteq P_{\nu}\subseteq RB_2^n$ and $P_{\nu}\to P$ in the Hausdorff metric, we have $\norm{\mathord{\cdot}}_{P_{\nu}}\to \norm{\mathord{\cdot}}_P$ locally uniformly on $\R^n$.

Fix $X\in C_c^{\infty}(\operatorname{int}P;\R^n)$, extended by zero to $\R^n$, and write $K_X=\operatorname{supp}X$. Thus, for sufficiently large $\nu$, the field $X$ is admissible in the dual energy formula for both $P_{\nu}$ and $P$. Recall that
\eq{
\label{eq:sec-3-hessian-passage}
\int_{\R^n}X^{\mathsf{T}}\,d(D^2\norm{\mathord{\cdot}}_{P_{\nu}})X
=\sum_{k,l=1}^n\int_{\R^n}
\norm{x}_{P_{\nu}}\partial_k\partial_l(X_kX_l)\,dx.
}
Applying \eqref{eq:sec-3-hessian-passage} with $P$ in place of $P_{\nu}$ and subtracting, we obtain
\eq{
&\biggl|\int_{\R^n}X^{\mathsf{T}}\,d(D^2\norm{\mathord{\cdot}}_{P_{\nu}})X
-\int_{\R^n}X^{\mathsf{T}}\,d(D^2\norm{\mathord{\cdot}}_P)X\biggr|\\
&\quad\leq \sup_{x\in K_X}\bigl|\norm{x}_{P_{\nu}}-\norm{x}_P\bigr|
\sum_{k,l=1}^n
\norm{\partial_k\partial_l(X_kX_l)}_{L^1(\R^n)}\to 0.
}
Moreover, the Cauchy--Schwarz inequality implies
\eq{
\biggl|\int_{\R^n}(\hat{f}_{\nu}-g)\operatorname{div}X\,dx\biggr|
\leq \norm{\hat{f}_{\nu}-g}_{L^2(B_R(0))}
\norm{\operatorname{div}X}_{L^2(\R^n)}\to 0.
}
Since $\mathsf{L}_{P_{\nu}}(f_{\nu},f_{\nu})$ converges along the chosen subsequence to the $\liminf$ of the original sequence, \eqref{eq:sec-3-dual-energy} and \eqref{eq:sec-3-radial-energy} imply
\eq{
&-2\int_Pg\operatorname{div}X\,dx
-\int_PX^{\mathsf{T}}\,d(D^2\norm{\mathord{\cdot}}_P)X\\
&\quad=\lim_{\nu\to \infty}
\left(-2\int_{P_{\nu}}\hat{f}_{\nu}\operatorname{div}X\,dx
-\int_{P_{\nu}}X^{\mathsf{T}}\,d(D^2\norm{\mathord{\cdot}}_{P_{\nu}})X\right)\\
&\quad\leq \frac{1}{n-1}\liminf_{\nu\to \infty}
\mathsf{L}_{P_{\nu}}(f_{\nu},f_{\nu}).
}
Taking the supremum over $X\in C_c^{\infty}(\operatorname{int}P;\R^n)$ and applying \eqref{eq:sec-3-dual-energy}, we obtain
\eq{
\label{eq:sec-3-dual-liminf}
\mathcal{L}_{\operatorname{int}P}(P,g)
\leq \frac{1}{n-1}\liminf_{\nu\to \infty}
\mathsf{L}_{P_{\nu}}(f_{\nu},f_{\nu})<\infty.
}

Let $(F_i)_{i\in J}$ be the facets of $P$. To show that $g$ is constant almost everywhere on the interior of each facet cone, fix $i\in J$ and consider
\eq{
W_i=\operatorname{int}(\operatorname{pos}F_i)\cap\operatorname{int}P.
}
The set $W_i$ is nonempty, open, and convex. Since $\norm{x}_P=\ip{v_i}{x}$ on $W_i$, the measure $D^2\norm{\mathord{\cdot}}_P$ vanishes on $W_i$. We test the dual energy formula with $tX$, where $X\in C_c^{\infty}(W_i;\R^n)$ and $t\in\R$, to obtain
\eq{
-2t\int_{W_i}g\operatorname{div}X\,dx
\leq \mathcal{L}_{\operatorname{int}P}(P,g)<\infty.
}
Thus
\eq{
\int_{W_i}g\operatorname{div}X\,dx=0
\quad\text{for every }X\in C_c^{\infty}(W_i;\R^n).
}
By \eqref{eq:sec-3-gradient-field-pairing}, $Dg=0$ on $W_i$. For a standard mollifier $\varphi_{\varepsilon}$ and any open ball $B_{\rho}(x)\Subset W_i$, we have $\nabla(\varphi_{\varepsilon}*g)=0$ on $B_{\rho}(x)$ for all small $\varepsilon> 0$. Passing to the $L^2$ limit as $\varepsilon\to 0^{+}$ shows that $g$ is constant almost everywhere on $B_{\rho}(x)$. Since $W_i$ is connected, $g=f_i$ almost everywhere on $W_i$ for a single constant $f_i$.

Now set $f=(f_i)_{i\in J}$. For every $m\in\mathbb{N}$, zero-homogeneity of $\hat{f}_{\nu}$ implies $g(m\mathord{\cdot})=g$ almost everywhere, and since $\operatorname{int}C_i=\bigcup_{m\in\mathbb{N}}mW_i$, we have $g=\hat{f}$ almost everywhere on $\R^n$. Thus \eqref{eq:sec-3-radial-L2-convergence} holds, and combining \eqref{eq:sec-3-dual-liminf} with \eqref{eq:sec-3-radial-energy} yields \eqref{eq:sec-3-energy-liminf-o}.

Finally, let $B=B_{R+1}(0)$. By \eqref{eq:sec-3-radial-L2-convergence}, we have $\hat{f}_{\nu}^k\to \hat{f}^k$ in $L^1(B)$ for $k=1,2$, while Hausdorff convergence ensures $\mathbf{1}_{P_{\nu}}\to \mathbf{1}_P$ almost everywhere. Since
\eq{
\bigl|\mathbf{1}_{P_{\nu}}\hat{f}_{\nu}^k-\mathbf{1}_P\hat{f}^k\bigr|
\leq \bigl|\hat{f}_{\nu}^k-\hat{f}^k\bigr|
+\bigl|\mathbf{1}_{P_{\nu}}-\mathbf{1}_P\bigr|\lvert\hat{f}^k\rvert,
}
applying the dominated convergence theorem to the second term yields $\mathbf{1}_{P_{\nu}}\hat{f}_{\nu}^k\to \mathbf{1}_P\hat{f}^k$ in $L^1(B)$. In view of \eqref{eq:sec-3-radial-mass}, the last two limits in \eqref{eq:sec-3-radial-L2-convergence} follow. If $P_{\nu}=-P_{\nu}$ and every $f_{\nu}$ is even, then $\hat{f}_{\nu}$ and $\hat{f}$ are even, and thus $f$ is even as well.
\end{proof}

\begin{proof}[Proof of \autoref{thm:polyhedral-lower-semicontinuity}]
There are constants $0<r<R$ such that
\eq{
rB_2^n\subseteq P_{\nu},\, P\subseteq RB_2^n
}
for all sufficiently large $\nu$. To prove lower semicontinuity, we may assume $r=1$ after a common dilation. If the $\liminf$ of the eigenvalues is infinite, then the conclusion is immediate. Otherwise, we may pass to a subsequence along which the eigenvalues converge to the $\liminf$ of the original sequence. For each $\nu$, choose a first nonconstant even eigenfunction $f_{\nu}$ satisfying
\eq{
\mathsf{M}_{P_{\nu}}(f_{\nu},\one)=0,
\quad
\mathsf{M}_{P_{\nu}}(f_{\nu},f_{\nu})=1,
\quad
\mathsf{L}_{P_{\nu}}(f_{\nu},f_{\nu})=\lambda_{1,e}(P_{\nu}).
}
The convergence of the eigenvalues along the subsequence implies that
\eq{
\mathsf{M}_{P_{\nu}}(f_{\nu},f_{\nu})
+\mathsf{L}_{P_{\nu}}(f_{\nu},f_{\nu})
=1+\lambda_{1,e}(P_{\nu})
}
is uniformly bounded. Therefore \autoref{prop:sec-3-facet-function-compactness} applies, and we can find an even facet function $f$ on $P$ such that
\eq{
\mathsf{M}_P(f,\one)=0,
\quad
\mathsf{M}_P(f,f)=1,
\quad
\mathsf{L}_P(f,f)\leq \liminf_{\nu\to \infty}\lambda_{1,e}(P_{\nu}).
}
In particular, $f$ is nonconstant and satisfies $\widetilde{\mathsf{M}}_P(f,f)=1$. Now \eqref{eq:sec-2-first-even-Rayleigh} yields
\eq{
\lambda_{1,e}(P)\leq \mathsf{L}_P(f,f)\leq \liminf_{\nu\to \infty}\lambda_{1,e}(P_{\nu}).
}

To prove upper semicontinuity, let us write $P=P(h)$ and choose a centered even eigenfunction $\tilde{f}$ for $\lambda_{1,e}(P)$ with $\mathsf{M}_P(\tilde{f},\tilde{f})=1$. Set $q_i=h_i(1+\varepsilon\tilde{f}_i)$ and $L=P(q)$. For sufficiently small $\varepsilon> 0$, each inequality $\ip{x}{u_i}\leq q_i$ in the representation of $P(q)$ defines a facet of $L$, and $h$ and $q$ lie in the closure of a common simple type chamber \cite[(8.6)]{HI26}. At each facet normal $u$ of $P_{\nu}$, define
\eq{
a_{\nu}(u)=\frac{h_L(u)}{h_{P_{\nu}}(u)}.
}
Note that the facet function $a_{\nu}$ is even. On the facet normals of $P$, the ratio $a=h_L/h_P$ satisfies $a_i=q_i/h_i=1+\varepsilon\tilde{f}_i$. The ratios $h_L/h_{P_{\nu}}$ converge uniformly to $h_L/h_P$ on $\Sn$, and the measures $h_{P_{\nu}}\,dS_{P_{\nu}}$ converge weakly to $h_P\,dS_P$. Since $\mathsf{M}_P(\tilde{f},\one)=0$, the formula for the centered mass form yields
\eq{
\widetilde{\mathsf{M}}_{P_{\nu}}(a_{\nu},a_{\nu})\to \widetilde{\mathsf{M}}_P(a,a)=\varepsilon^2\mathsf{M}_P(\tilde{f},\tilde{f})=\varepsilon^2.
}

By \cite[Lem. 8.9]{HI26} and continuity of mixed volumes,
\eq{
\limsup_{\nu\to \infty}\mathsf{L}_{P_{\nu}}(a_{\nu},a_{\nu})
&\leq (n-1)\int_{\Sn}\frac{h_L^2}{h_P}\,dS_P-n(n-1)\mathcal{V}(P[n-2],L[2])\\
&=\mathsf{L}_P(a,a)=\varepsilon^2\lambda_{1,e}(P),
}
where $\mathcal{V}(P[n-2],L[2])$ denotes the mixed volume of $n-2$ copies of $P$ and two copies of $L$. To verify the first equality, let $V(h)=\Vol(P(h))$. Since $h$ and $q$ belong to a common chamber closure, \cite[(8.28)]{HI26} yields $D^2V(h)[q,q]=n(n-1)\mathcal{V}(P[n-2],L[2])$. Since $q_i=h_i a_i$, \cite[(8.11)]{HI26} and the definition of the facet masses yield
\eq{
\mathsf{L}_P(a,a)=(n-1)\mathsf{M}_P(a,a)-D^2V(h)[q,q],\quad
\mathsf{M}_P(a,a)=\int_{\Sn}\frac{h_L^2}{h_P}\,dS_P.
}
Finally, $a=\one+\varepsilon\tilde{f}$ and constants have zero energy.

Since the facet function $a_{\nu}$ is nonconstant for all sufficiently large $\nu$, by \eqref{eq:sec-2-first-even-Rayleigh},
\eq{
\limsup_{\nu\to \infty}\lambda_{1,e}(P_{\nu})
\leq \limsup_{\nu\to \infty}\frac{\mathsf{L}_{P_{\nu}}(a_{\nu},a_{\nu})}{\widetilde{\mathsf{M}}_{P_{\nu}}(a_{\nu},a_{\nu})}
\leq \lambda_{1,e}(P).
}
\end{proof}

\begin{corollary}\label{cor:radial-eigenfunction-convergence}
Under the hypotheses of \autoref{thm:polyhedral-lower-semicontinuity}, let $f_{\nu}$ be a centered unit-mass even eigenfunction for $\lambda_{1,e}(P_{\nu})$. Then there is a subsequence along which
\eq{
\hat{f}_{\nu}\to \hat{f}\quad\text{in }L^2_{\mathrm{loc}}(\R^n)
}
for some centered unit-mass even eigenfunction $f$ of $P$ with eigenvalue $\lambda_{1,e}(P)$.
\end{corollary}

\begin{proof}
By \autoref{prop:sec-3-facet-function-compactness}, there is a subsequence along which $\hat{f}_{\nu}\to \hat{f}$ in $L^2_{\mathrm{loc}}(\R^n)$ for a centered even facet function $f$ of unit mass. Due to \eqref{eq:sec-2-first-even-Rayleigh} and \autoref{thm:polyhedral-lower-semicontinuity},
\eq{
\lambda_{1,e}(P)\leq \mathsf{L}_P(f,f)\leq \lim_{\nu\to \infty}\lambda_{1,e}(P_{\nu})=\lambda_{1,e}(P).
}
Thus $f$ is an eigenfunction for $\lambda_{1,e}(P)$.
\end{proof}

\section{Existence of a minimizer under a facet bound}

We first show that $\lambda_{1,e}(P)$ is invariant under invertible linear maps. For the corresponding result for smooth convex bodies, see \cite[Thm. 5.8]{KM22}.

\begin{lemma}\label{lem:linear-invariance}
Let $P\subseteq\R^n$ be a full-dimensional origin-symmetric polytope. If $\ell:\R^n\to \R^n$ is an invertible linear map, then $\lambda_{1,e}(\ell P)=\lambda_{1,e}(P)$.
\end{lemma}

\begin{proof}
Let us index the facets of $P$ and $\ell P$ as $(F_i)_{i\in I}$ and $(\ell F_i)_{i\in I}$ by a signed label set $I$. To each facet function $f=(f_i)_{i\in I}$ on $P$, associate the facet function $g=(g_i)_{i\in I}$ on $\ell P$ defined by $g_i=f_i$ for every $i\in I$. Since $\ell$ maps each facet cone $\operatorname{pos}F_i$ onto $\operatorname{pos}(\ell F_i)$, the radial extensions satisfy $\hat{g}(y)=\hat{f}(\ell^{-1}y)$. 

To compare the dual energies of $\hat{f}$ and $\hat{g}$, set $q(x)=\norm{x}_P$ and $\widetilde{q}(y)=\norm{y}_{\ell P}=q(\ell^{-1}y)$. To each $X\in C_c^{\infty}(\operatorname{int}P;\R^n)$, we associate the vector field $Y\in C_c^{\infty}(\operatorname{int}(\ell P);\R^n)$ defined by $Y(y)=\ell X(\ell^{-1}y)$. Note that $\operatorname{div}Y(\ell x)=\operatorname{div}X(x)$, and after the change of variables $y=\ell x$:
\eq{
\int_{\operatorname{int}(\ell P)}\hat{g}(y)\operatorname{div}Y(y)\,dy
=|\det\ell|\int_{\operatorname{int}P}\hat{f}(x)\operatorname{div}X(x)\,dx.
}

Let $(\ell_{ka})$ and $((\ell^{-1})_{rk})$ denote the matrices of $\ell$ and $\ell^{-1}$ in the standard basis of $\R^n$. For any $\psi\in C_c^{\infty}(\R^n)$,
\eq{
(\partial_{y_k}\psi)(\ell x)
=\sum_{r=1}^n(\ell^{-1})_{rk}\partial_{x_r}\bigl(\psi(\ell x)\bigr).
}
Now, since $Y_k(\ell x)=\sum_{a=1}^n\ell_{ka}X_a(x)$, applying the chain rule twice to $Y_kY_l$ and summing over $k$ and $l$, we obtain
\eq{
\sum_{k,l=1}^n\bigl(\partial_{y_k}\partial_{y_l}(Y_kY_l)\bigr)(\ell x)
&=\sum_{k,l,r,s,a,b=1}^n(\ell^{-1})_{rk}\ell_{ka}(\ell^{-1})_{sl}\ell_{lb}\partial_{x_r}\partial_{x_s}\bigl(X_a(x)X_b(x)\bigr)\\
&=\sum_{a,b=1}^n\partial_{x_a}\partial_{x_b}\bigl(X_a(x)X_b(x)\bigr).
}
Therefore, by \eqref{eq:sec-3-distributional-definitions} and the change of variables $y=\ell x$,
\eq{
\int_{\operatorname{int}(\ell P)}Y^{\mathsf{T}}\,d(D^2\widetilde{q})Y
&=\sum_{k,l=1}^n\int_{\R^n}\widetilde{q}(y)\partial_{y_k}\partial_{y_l}\bigl(Y_k(y)Y_l(y)\bigr)\,dy\\
&=|\det\ell|\int_{\R^n}q(x)\sum_{k,l=1}^n\bigl(\partial_{y_k}\partial_{y_l}(Y_kY_l)\bigr)(\ell x)\,dx\\
&=|\det\ell|\sum_{a,b=1}^n\int_{\R^n}q(x)\partial_{x_a}\partial_{x_b}\bigl(X_a(x)X_b(x)\bigr)\,dx=|\det\ell|\int_{\operatorname{int}P}X^{\mathsf{T}}\,d(D^2q)X.
}

Since $X\mapsto Y$ is a bijection between $C_c^{\infty}(\operatorname{int}P;\R^n)$ and $C_c^{\infty}(\operatorname{int}(\ell P);\R^n)$, it follows from \eqref{eq:sec-3-dual-energy} that
\eq{
\mathcal{L}_{\operatorname{int}(\ell P)}(\ell P,\hat{g})
=|\det\ell|\mathcal{L}_{\operatorname{int}P}(P,\hat{f}).
}
Similarly, changing variables in \eqref{eq:sec-3-radial-mass} yields
\eq{
\mathsf{M}_{\ell P}(g,g)=|\det\ell|\mathsf{M}_P(f,f),\quad
\mathsf{M}_{\ell P}(g,\one)=|\det\ell|\mathsf{M}_P(f,\one).
}
Finally, since $f\mapsto g$ is a bijection preserving evenness and nonconstancy, \eqref{eq:sec-3-radial-energy} and \eqref{eq:sec-2-first-even-Rayleigh} imply $\lambda_{1,e}(\ell P)=\lambda_{1,e}(P)$.
\end{proof}

For an integer $N\geq n$, let $\mathcal{P}_N$ be the class of full-dimensional origin-symmetric polytopes with at most $N$ pairs of facets, and define
\eq{
\label{eq:bounded-facet-infimum}
\lambda_N^{\star}=\inf_{P\in\mathcal{P}_N}\lambda_{1,e}(P).
}
The class $\mathcal{P}_N$ is nonempty because it contains every origin-symmetric parallelotope.

\begin{corollary}\label{cor:polyhedral-spectral-gap}
Every full-dimensional origin-symmetric polytope $P\subseteq\R^n$ satisfies
\eq{
\lambda_{1,e}(P)\geq n-1+\frac{1}{n-1}.
}
\end{corollary}

\begin{proof}
Let $N$ be the number of pairs of facets of $P$, and choose $(P_{\nu})\subseteq\mathcal{P}_N$ such that $\lambda_{1,e}(P_{\nu})\to \lambda_N^{\star}$. By John's theorem \cite[p. 588, Thm. 10.12.2]{Sch14}, for each $\nu$ there is an invertible linear map $\ell_{\nu}$ such that $B_2^n\subseteq\ell_{\nu}P_{\nu}\subseteq\sqrt{n}B_2^n$. Since invertible linear maps preserve the number of facet pairs, in view of \autoref{lem:linear-invariance}, we may replace $P_{\nu}$ by $\ell_{\nu}P_{\nu}$. Denote the normalized sequence again by $(P_{\nu})$. By Blaschke's selection theorem \cite[p. 63, Thm. 1.8.7]{Sch14}, after passing to a subsequence, $P_{\nu}\to P_{\star}$ in the Hausdorff metric. The ball inclusions and origin symmetry pass to $P_{\star}$. Note that $P_{\nu}^{\ast}\to P_{\star}^{\ast}$ in the Hausdorff metric and each $P_{\nu}^{\ast}\subseteq B_2^n$ is the convex hull of $N$ (possibly repeated) antipodal pairs of vertices. Passing to a further subsequence along which each pair converges shows that $P_{\star}^{\ast}$ is a polytope with at most $N$ pairs of vertices, and thus $P_{\star}=(P_{\star}^{\ast})^{\ast}\in\mathcal{P}_N$. Therefore, by \autoref{thm:polyhedral-lower-semicontinuity} and \eqref{eq:bounded-facet-infimum},
\eq{
\lambda_N^{\star}
\leq \lambda_{1,e}(P_{\star})
\leq \liminf_{\nu\to \infty}\lambda_{1,e}(P_{\nu})
=\lambda_N^{\star}.
}
In a complete irredundant facet representation of $P_{\star}$ with signed label set $I$, write $P_{\star}=P(h_{\star})$ with even $h_{\star}\in\mathcal{U}_{\mathrm{irr}}^I$. By \cite[Lem. 8.2]{HI26}, $\mathcal{U}_{\mathrm{irr}}^I$ is open in $\R^I$. Moreover, since $|I|\leq 2N$, every even $q\in\mathcal{U}_{\mathrm{irr}}^I$ satisfies $P(q)\in\mathcal{P}_N$ and $\lambda_{1,e}(P(q))\geq \lambda_N^{\star}=\lambda_{1,e}(P_{\star})$. Thus $h_{\star}$ is a local minimizer of $\lambda_{1,e}$ under even support variations, and \cite[Thm. 1.1]{HI26} yields
\eq{
\lambda_{1,e}(P)\geq \lambda_N^{\star}=\lambda_{1,e}(P_{\star})\geq n-1+\frac{1}{n-1}.
}
\end{proof}

\section{\texorpdfstring{The $L^{\frac{n-2}{n-1}}$-Brunn--Minkowski inequality}{The Lp-Brunn--Minkowski inequality at p=(n-2)/(n-1)}}

Set $p_n=1-1/(n-1)$ and $c_n=n-p_n=n-1+1/(n-1)$. For the next theorem, we recall the mixed-volume expansion \cite[p. 280, Thm. 5.1.7]{Sch14}
\eq{
\Vol(K+tL)=\sum_{k=0}^n\binom{n}{k}\mathcal{V}(K[n-k],L[k])t^k,
\quad t\geq 0.
}
For $1\leq k\leq n-1$, the bracketed multiplicities in the mixed-volume notation mean
\eq{
\mathcal{V}(K[n-k],L[k])
=\mathcal{V}\bigl(\underbrace{K,\ldots,K}_{n-k\text{ times}},
\underbrace{L,\ldots,L}_{k\text{ times}}\bigr).
}
When $k=1$, we abbreviate $\mathcal{V}(K[n-1],L[1])$ as $\mathcal{V}(K[n-1],L)$, and $\mathcal{V}(K[n])=\Vol(K)$. Mixed volumes are continuous under Hausdorff convergence \cite[pp. 280--281]{Sch14}.

\begin{theorem}\label{thm:local-Lpn-inequality}
For all origin-symmetric convex bodies $K,L\subseteq\R^n$,
\eq{
\label{eq:local-Lpn-inequality}
\Vol(K)\left(n(n-1)\mathcal{V}(K[n-2],L[2])+\frac{1}{n-1}\int_{\Sn}\frac{h_L^2}{h_K}\,dS_K\right)\leq nc_n\mathcal{V}(K[n-1],L)^2.
}
\end{theorem}

\begin{proof}
We first prove the inequality when $K=P$ is a polytope. In a complete irredundant facet representation of $P$ with signed label set $I$, put $a_i=h_L(u_i)/h_i$. Since $h_L(-u_i)=h_L(u_i)$ and $h_{-i}=h_i$, the facet function $a$ is even. Using the first mixed-volume formula \cite[p. 280, Thm. 5.1.7]{Sch14}, we obtain
\eq{
\mathsf{M}_P(\one,\one)
&=\sum_{i\in I}h_i\mathcal{H}^{n-1}(F_i)
=\int_{\Sn}h_P\,dS_P=n\Vol(P),\\
\mathsf{M}_P(a,\one)
&=\sum_{i\in I}h_L(u_i)\mathcal{H}^{n-1}(F_i)
=n\mathcal{V}(P[n-1],L),\\
\mathsf{M}_P(a,a)
&=\sum_{i\in I}\frac{h_L(u_i)^2}{h_i}\mathcal{H}^{n-1}(F_i)
=\int_{\Sn}\frac{h_L^2}{h_P}\,dS_P.
}
The centered function $a-\bar{a}\one$ is even and satisfies
\eq{
\mathsf{M}_P(a-\bar{a}\one,\one)&=0,\\
\mathsf{M}_P(a-\bar{a}\one,a-\bar{a}\one)
&=\mathsf{M}_P(a,a)-\frac{\mathsf{M}_P(a,\one)^2}{\mathsf{M}_P(\one,\one)},\\
\mathsf{L}_P(a-\bar{a}\one,a-\bar{a}\one)&=\mathsf{L}_P(a,a).
}
If $a$ is nonconstant, we apply \eqref{eq:sec-2-first-even-Rayleigh} and \autoref{cor:polyhedral-spectral-gap} to $a-\bar{a}\one$. If $a$ is constant, then $a-\bar{a}\one=0$ and $\mathsf{L}_P(a,a)=0$. In either case, it follows that
\eq{
c_n\left(\mathsf{M}_P(a,a)-\frac{\mathsf{M}_P(a,\one)^2}{\mathsf{M}_P(\one,\one)}\right)
\leq \mathsf{L}_P(a,a).
}
In view of \cite[Lem. 8.9]{HI26},
\eq{
&c_n\left(\int_{\Sn}\frac{h_L^2}{h_P}\,dS_P
-\frac{n\mathcal{V}(P[n-1],L)^2}{\Vol(P)}\right)\\
&\leq \mathsf{L}_P(a,a)\leq (n-1)\int_{\Sn}\frac{h_L^2}{h_P}\,dS_P-n(n-1)\mathcal{V}(P[n-2],L[2]).
}
Rearranging, we conclude that
\eq{
\Vol(P)\left(n(n-1)\mathcal{V}(P[n-2],L[2])
+\frac{1}{n-1}\int_{\Sn}\frac{h_L^2}{h_P}\,dS_P\right)\leq nc_n\mathcal{V}(P[n-1],L)^2.
}

For general $K$, choose full-dimensional origin-symmetric polytopes $P_{\nu}$ converging to $K$ in the Hausdorff metric, and keep $L$ fixed. Passing to the limit in the inequality for $P_{\nu}$, we obtain \eqref{eq:local-Lpn-inequality}.
\end{proof}

\begin{definition}
For $0<p<1$ and origin-symmetric convex bodies $K,L\subseteq\R^n$, define
\eq{
(1-t)K+_ptL=\bigcap_{u\in\Sn}\left\{x\in\R^n:\ip{x}{u}\leq \bigl((1-t)h_K(u)^p+th_L(u)^p\bigr)^{\frac{1}{p}}\right\},\quad 0\leq t\leq 1.
}
\end{definition}

\begin{corollary}\label{cor:global-Lpn-inequality}
Let $p_n=1-1/(n-1)$. For all origin-symmetric convex bodies $K,L\subseteq\R^n$ and every $0\leq t\leq 1$,
\eq{
\Vol\bigl((1-t)K+_{p_n}tL\bigr)^{\frac{p_n}{n}}\geq (1-t)\Vol(K)^{\frac{p_n}{n}}+t\Vol(L)^{\frac{p_n}{n}}.
}
\end{corollary}

\begin{proof}
Note that $n-p_n=c_n$. \autoref{thm:local-Lpn-inequality} asserts that
\eq{
&\Vol(K)\left(n(n-1)\mathcal{V}(K[n-2],L[2])
+(1-p_n)\int_{\Sn}\frac{h_L^2}{h_K}\,dS_K\right)\\
&\quad-n(n-p_n)\mathcal{V}(K[n-1],L)^2\leq 0
}
for every pair of origin-symmetric convex bodies $K,L\subseteq\R^n$. Thus \eqref{eq:local-Lpn-inequality} is precisely the statement of \cite[Conj. 1.6]{Put21} at the exponent $p=p_n$. By \cite[Thm. 1.7]{Put21}, the global $L^{p_n}$-Brunn--Minkowski inequality holds.
\end{proof}

\begin{remark}
For an origin-symmetric $C^2_{+}$ convex body $K\subseteq\R^n$, we recall from \cite[Sec. 1]{HI26} the centro-affine metric and volume measure
\eq{
g_K=\frac{r_K}{h_K},\quad dV_K=h_K\det\nolimits_{\bar{g}}(r_K)\,d\omega,
}
where the determinant is computed in a $\bar{g}$-orthonormal frame. The first nonconstant even eigenvalue of $K$ has the Rayleigh characterization
\eq{
\lambda_{1,e}(K)=\inf_{\substack{0\neq \psi\in C^2(\Sn),\ \psi(-u)=\psi(u)\\
\int_{\Sn}\psi\,dV_K=0}}
\frac{\displaystyle\int_{\Sn}|\nabla\psi|_{g_K}^2\,dV_K}
{\displaystyle\int_{\Sn}\psi^2\,dV_K}.
}
The gradient is taken with respect to $g_K$, with
\eq{
|\nabla\psi|_{g_K}^2=g_K^{-1}(d\psi,d\psi)=h_Kr_K^{-1}(d\psi,d\psi).
}
By polyhedral approximation \cite[p. 67, Thm. 1.8.16]{Sch14}, \autoref{cor:polyhedral-spectral-gap}, and upper semicontinuity \cite[Thm. 8.11]{HI26}, every origin-symmetric $C^2_{+}$ convex body $K\subseteq\R^n$ satisfies $\lambda_{1,e}(K)\geq c_n$.
\end{remark}

The logarithmic Brunn--Minkowski inequality, corresponding to $p=0$, was proved for origin-symmetric planar bodies by B\"or\"oczky, Lutwak, Yang, and Zhang \cite[Thm. 1.3]{BLYZ12}, and for pairs of bodies that are unconditional with respect to a common orthonormal basis by Saroglou \cite[Thm. 1.2]{Sar15}. The planar and unconditional logarithmic Brunn--Minkowski inequalities also imply the global $L^p$-Brunn--Minkowski inequality for every $0<p<1$ in the respective classes; see \cite[Thm. 1.7]{BLYZ12} and \cite[Cor. 1.3]{Sar15}. In arbitrary dimension, the local $L^p$-Brunn--Minkowski inequality for smooth origin-symmetric convex bodies was established in \cite[Thm. 1.1]{KM22} for $1-c_0n^{-\frac{3}{2}}\leq p<1$, where $c_0> 0$ is universal. The local-to-global results in \cite[Thm. 1.1]{CHLL20} and \cite[Thm. 1.7, Cor. 1.8]{Put21} imply the global inequality in the same range.

B\"or\"oczky and Kalantzopoulos proved the logarithmic Brunn--Minkowski inequality for convex bodies invariant under the same $n$ linear reflections whose fixed hyperplanes intersect only at the origin \cite[Thm. 2]{BK22}. The proof reduces the problem to the unconditional case through Weyl chambers. B\"or\"oczky and De established stability estimates for the logarithmic Brunn--Minkowski and logarithmic Minkowski inequalities under the same symmetry assumptions \cite[Thms. 1.4, 1.5]{BD24}. More recently, Lu proved the global logarithmic Brunn--Minkowski inequality for origin-symmetric bodies invariant under the same $n-2$ hyperplane reflections with pairwise orthogonal normal vectors \cite[Thm. 1.4]{Lu26}. Also, \cite[Cor. 1.5]{Lu26} covers pairs of origin-symmetric bodies of revolution with arbitrary axes.

For an origin-symmetric zonoid $K$, van Handel proved the local logarithmic Brunn--Minkowski inequality and the logarithmic Minkowski inequality, both with an arbitrary origin-symmetric convex body $L$ \cite[Thm. 1.4, Cor. 1.6]{vH23}. Iffland proved the local logarithmic Brunn--Minkowski inequality when $K$ is an origin-symmetric body of revolution and $L$ is an arbitrary origin-symmetric convex body \cite[Thm. A]{Iff26}. For unconditional $C^{\infty}_{+}$ convex bodies $K$, the centro-affine Poincar\'e inequality implies $\lambda_{1,e}(K)\geq n$ \cite[Thm. 1.1]{HIUnc26}.

The even $L^p$-Minkowski inequality and uniqueness in the even $L^p$-Minkowski problem were established for a range of $p$ determined by curvature pinching assumptions \cite[Thm. 1.2]{IM24}. Milman proved an isomorphic form of the logarithmic Minkowski inequality: every origin-symmetric convex body $K$ admits an origin-symmetric $C^{\infty}_{+}$ convex body $L$ with $K\subseteq L\subseteq 8K$ such that $L$ satisfies the logarithmic Minkowski inequality against every origin-symmetric convex body \cite[Cor. 1.6]{Mil25}. The cone-volume measure of $L$ determines $L$ uniquely among origin-symmetric convex bodies.

\section*{Acknowledgments}
Hu was supported by the National Key Research and Development Program of China (Grant No. 2021YFA1001800). Ivaki was supported by the Austrian Science Fund (FWF) under Project P36545.

\section*{AI disclosure}
The authors used AI tools during the development of this paper.

\bigskip
\begingroup
\footnotesize
\noindent
\begin{minipage}[t]{0.40\textwidth}
\textsc{School of Mathematical\\ Sciences, Beihang University,\\ Beijing 100191, China}\\
\href{mailto:huyingxiang@buaa.edu.cn}{\texttt{huyingxiang@buaa.edu.cn}}
\end{minipage}
\hfill
\begin{minipage}[t]{0.40\textwidth}
\textsc{Institut f\"{u}r Diskrete\\ Mathematik und Geometrie,\\ Technische Universit\"{a}t Wien,\\ Wiedner Hauptstra{\ss}e 8--10,\\ 1040 Wien, Austria}\\
\href{mailto:mohammad.ivaki@tuwien.ac.at}{\texttt{mohammad.ivaki@tuwien.ac.at}}
\end{minipage}
\par
\endgroup

\end{document}